\documentclass[10pt]{article}

\usepackage{amsmath}
\usepackage{amsfonts}
\usepackage{amsthm}
\usepackage[english]{babel}
\usepackage{graphicx}
\usepackage[all]{xy}
\usepackage{amssymb}
\usepackage{bm}

\newtheorem{theorem}{Theorem}[section]
\newtheorem{definition}[theorem]{Definition}

\newtheorem{remark}[theorem]{Remark}

\newtheorem{notation}[theorem]{Notation}

\newcommand{\mR}{\mathbb{R}}

\newcommand{\mH}{\mathbb{H}}

\newcommand{\cF}{\mathcal{F}}

\newcommand{\cI}{\mathcal{I}}

\newcommand{\qi}{\bm{i}}
\newcommand{\qj}{\bm{j}}
\newcommand{\qk}{\bm{k}}

\begin{document}

\title{Connecting spatial and frequency domains for the quaternion Fourier transform}

\author{H.\ De Bie\footnote{Corresponding author; E-mail: {\tt Hendrik.DeBie@UGent.be}} \and N.\ De Schepper\footnote{E-mail: {\tt nds@cage.ugent.be}} \and T.A. Ell\footnote{E-mail: {\tt t.ell@ieee.org}}\and  K. Rubrecht\footnote{E-mail: {\tt klausrubrecht@gmail.com}} \and S.J. Sangwine\footnote{E-mail: {\tt sjs@essex.ac.uk}}}

\vspace{10mm}
\date{\small{H. De Bie and K. Rubrecht:  Department of Mathematical Analysis, Faculty of Engineering and Architecture\\}
\small{Ghent University, Galglaan 2, 9000 Gent,
Belgium}\\ \vspace{5mm}
\small{N. De Schepper: Department of Mathematics, University of Antwerp\\ Middelheimlaan 2, 2020 Antwerp, Belgium.}\\
\vspace{5mm}
\small{T.A. Ell: UTC Aerospace Systems, 14300 Judicial Rd, Burnsville, MN 55306-4898, USA.}\\
\vspace{5mm}
\small{S.J. Sangwine: School of Computer Science and Electronic Engineering, University of Essex\\ Wivenhoe Park, Colchester CO4 3SQ, United Kingdom}
}

\maketitle

\begin{abstract}
The quaternion Fourier transform (qFT) is an important tool in multi-dimensional data analysis, in particular for the study of color images. An important problem when applying the qFT is the mismatch between the spatial and frequency domains: the convolution of two quaternion signals does not map to the pointwise product of their qFT images. The recently defined `Mustard'
convolution behaves nicely in the frequency domain, but complicates the 
corresponding spatial domain analysis. 

The present paper analyses in detail the correspondence between classical convolution and the new Mustard convolution. In particular, an expression is derived that allows one to write classical convolution as a finite linear combination of suitable Mustard convolutions. This result is expected to play a major role in the further development of quaternion image processing, as it yields a formula for the qFT spectrum of the classical convolution.
\end{abstract}

\noindent
\textbf{Keywords:} Quaternion Fourier transform; convolution products; frequency domain; spatial domain

\section{Introduction}
\label{sec:1}

Since more than a decade ago, the use of quaternions has been investigated for 3D data analysis. As a guiding example we use color images. Indeed, one can represent, say, a color image as a pure quaternion via the identification
\[
R(x,y) \qi + G(x,y) \qj + B(x,y) \qk
\]
where $x$ and $y$ denote the pixel coordinates, and $R(x,y)$, $G(x,y)$ and $B(x,y)$ are the red, green and blue color channels. When pixels are represented in such a way, they can subsequently be processed in a non-marginal way. This means that all color components of the image are treated together, as opposed to processing each of the three components independently. The main mathematical tool applied in this context is the quaternion Fourier transform (qFT) \cite{EPHD, S, BS, PDC, ES}. This transform generalizes the classical Fourier transform by replacing the imaginary unit by two pure quaternions $\mu$ and $\nu$ that square to -1. These two quaternions can be chosen independently, but are typically related to relevant directions in the color space, such as the gray line or a direction perpendicular to it.

So far, the quaternion approach has led to many successful applications. We mention edge detection and image filtering \cite{S_CE, DCF}, flow visualization \cite{Ebl2}, watermarking \cite{BLBC}, pattern recognition \cite{MST, SC}, generalized wavelets \cite{BC, BC2, BTN}, noise removal from video images \cite{J} and many more. In particular for saliency detection, which is the computerised detection of image features likely to attract human attention, use of a qFT method based on its phase spectrum vastly outperforms other methods in terms of efficiency and robustness \cite{GZ2}.

The success of classical complex signal analysis is almost entirely due to the interplay between the Fourier transform (FT) and the classical convolution product and its application to the study of linear time-invariant (LTI) systems. Design of filters is very flexible and can be done both in the time and frequency domains. A fast computer implementation is obtained in the frequency domain, as the FT maps the computationally involved convolution of two signals to the much simpler pointwise product of their respective Fourier representations.

The main issue in applying the qFT, as already indicated in \cite{PDC}, is the mismatch between spatial and frequency domain: the qFT of the convolution of two quaternion signals does not equal the product of their respective qFTs. This means that, although e.g. the qFT saliency detection method works extremely well, it is not understood why it works so well. In turn, this means that other applications are not easily derived. The situation is even worse for the quaternion color edge filter \cite{S_CE}, which, by exploiting the geometric properties of the quaternions, generalizes the Prewitt edge detector and detects sharp changes of color in images. This filter is still implemented using the classical convolution product, thus leading to a slow algorithm, as the qFT cannot yet be applied in a meaningful way.

This brings us to the aim of the present paper, namely to correctly connect spatial and frequency domain for the qFT, i.e. we want to match quaternion filter design in the spatial domain (using classical convolution) with design in the quaternion frequency domain (using pointwise products). This will allow us, in a series of follow-up papers, to construct a complete signal processing theory for quaternion LTI systems, which would be a huge breakthrough. 

The key idea to solve this problem is to find a suitable decomposition of the two signals to be convolved, so that the action of the qFT simplifies to a sum of products. This will be tackled using a mathematical tool we call Mustard convolution, introduced in \cite{BDDS} as follows:
\[
 f*_{\mu,\nu}g:=2\pi (\cF^{\mu,\nu})^{-1}(\cF^{\mu,\nu}(f)\cF^{\mu,\nu}(g)).
\]
with $f$ and $g$ two quaternion signals and $\cF^{\mu,\nu}$ the qFT. We will establish a precise connection between classical convolution and Mustard convolution by expressing the standard convolution $\ast$ as a finite linear combination of  Mustard convolutions $*_{\mu,\nu}$ in the following way:
\begin{equation}
\label{MCCC}
f \ast g = \sum_i (a_i f^{\phi_i}) *_{\mu,\nu} (b_i g^{\psi_i}), \qquad a_i, b_i \in \mH.
\end{equation}
where $\phi_i, \psi_i$ denote suitable reflections on the arguments of the signals $f$ and $g$ (see Notation \ref{not:phi,gamma}). Subsequently applying the qFT to this formula then yields the desired connection between spatial domain and frequency domain, see e.g. Theorem \ref{CCspec} for the case of perpendicular roots.

The expression (\ref{MCCC})  will be derived in three subsequent steps, depending on the relation between $\mu$ and $\nu$, the quaternion roots of -1 used to construct the qFT. We will distinguish
\begin{itemize}
\item the case where $\mu = \nu$, see Theorem \ref{th:convolutiontheorem1}
\item the case where $\mu$ is perpendicular to $\nu$, see Theorem \ref{th:convolutiontheorem2}
\item the general case, see Theorem \ref{th:convolutiontheorem3}.
\end{itemize}

Whereas the case $\mu = \nu$ is relatively simple to deal with, the proof of the general case is hard. This lies in the fact that the quaternions do not form a commutative algebra and extra care has to be taken in order to swap exponentials, the roots  $\mu$ and $\nu$, and the quaternion signals $f$ and $g$ in order to arrive at the desired result. 
As the resulting theorems are complicated for arbitrary choices of roots, we have implemented and thoroughly verified all our results in Theorems \ref{th:mustclassconv}, \ref{th:convolutiontheorem1}, \ref{th:convolutiontheorem2} and \ref{th:convolutiontheorem3} using the symbolic algebra software Maple.

The paper is organised as follows. In Section \ref{sec:2} we give basic definitions for quaternions and the qFT and derive several important identities we will need later on. In Section \ref{sec:3} we establish our result for equal and perpendicular roots. Next, in Section \ref{sec:gen} we treat the general case of a qFT defined using two arbitrary roots of -1. Finally, in Section \ref{sec:con} we discuss some consequences of our results.

\section{Quaternions and the quaternion Fourier transform}
\label{sec:2}
\subsection{Quaternions}
The algebra of the quaternions $\mH$ is defined as follows:

\begin{definition}
The quaternion algebra $\mH$  is the algebra over $\mR$ generated by the basis elements $\qi,\qj,\qk$, under the relations
\begin{equation*}
\label{eq:basicquaternions}
\qi^2=\qj^2=\qk^2=\qi\qj\qk=-1.
\end{equation*}
These relations immediately imply: 
\begin{displaymath}
\qi\qj = -\qj\qi = \qk, \qquad \qj\qk = -\qk\qj = \qi, \qquad \qk\qi = -\qi\qk = \qj.
\end{displaymath}
\end{definition}

A general quaternion $q$ is characterized by four real numbers:
\begin{equation*}
q= a + b\qi +c\qj + d\qk. 
\end{equation*}
Sometimes it is useful to split a quaternion into a scalar and a vector part. The scalar part is the real part of the quaternion, $S(q)= a$. The vector part of the quaternion consists of the purely quaternionic part $V(q) = b\qi + c\qj + d\qk$. The quaternion product of $q_1=  a_1 + b_1\qi +c_1\qj + d_1\qk$ and $q_2=  a_2 + b_2\qi +c_2\qj + d_2\qk$ is given by  
\begin{equation*}
\begin{aligned}
q_1 q_2 &= a_1a_2 - b_1b_2 - c_1c_2 - d_1d_2 +(a_1b_2 + b_1a_2 + c_1d_2 - d_1c_2)\qi +(a_1c_2 - b_1d_2 + c_1a_2 + d_1b_2)\qj\\
&+(a_1d_2 + b_1c_2 - c_1b_2 + d_1a_2)\qk.
\end{aligned}
\end{equation*}

In the rest of the paper, we will always consider quaternion functions. These are functions $f$ from $\mR^2$ to $\mH$. They can be expanded as
\[
f = f_0 + f_1 \qi +f_2 \qj + f_3 \qk,
\]
where we will assume $f_i \in L^1(\mR^2)$ so that Fourier transform and convolution are always well-defined.

\subsection{The quaternion Fourier transform}
In this paper we consider the left quaternion Fourier transform which is defined as follows:

\begin{definition}\label{def:qft}
Denote by $\mu,\nu \in \cI_2=\lbrace a \in \mH | a^2=-1 \rbrace$ two geometric square roots of minus one, i.e. $\mu^2=\nu^2=-1$. The left qFT of a function $f: \mR^{2} \rightarrow  \mH$ is defined as:
\begin{equation}
\label{eq:Transform}
\mathcal{F}^{\mu,\nu} (f) (y) = \frac{1}{2\pi} \int_{\mathbb{R}^2} e^{-\mu x_1 y_1} e^{-\nu x_2 y_2} f(x) dx
\end{equation}
with $x = (x_1,x_2)$ to simplify notations.
\end{definition}

The inverse transform is given by
\[
\left( \mathcal{F}^{\mu,\nu} \right)^{-1} (f) (y) = \frac{1}{2\pi} \int_{\mathbb{R}^2}  e^{\nu x_2 y_2} e^{\mu x_1 y_1} f(x) dx.
\]

There exist a plethora of different Fourier transforms in quaternion analysis, and the above definition is obviously not unique. Due to the non-commutative  character of the quaternions the transform can also be defined by sandwiching the function $f$ with two exponentials, or by putting the exponentials on the right side. These are not the only three options, for an overview of different qFTs and their basic properties the reader can for example consult  \cite{BR1, E2}. In the present paper, we will always work with the left qFT, but the computations can be adapted to other definitions without additional difficulties.

Based on an idea of Mustard in the context of the fractional Fourier transform \cite{Must}, we define a generalized convolution for the qFT following \cite{BDDS}. 
\begin{definition}\label{def:conv_must}
For the qFT $\cF^{\mu,\nu}$, the Mustard convolution $*_{\mu,\nu}$ is given by
\begin{equation*}
 (f*_{\mu,\nu}g)(x):=2\pi (\cF^{\mu,\nu})^{-1}(\cF^{\mu,\nu}(f)\cF^{\mu,\nu}(g))(x).
\end{equation*}
\end{definition}
It clearly satisfies
\[
\cF^{\mu,\nu} (f*_{\mu,\nu}g) = 2\pi\cF^{\mu,\nu}(f)\cF^{\mu,\nu}(g).
\]

We introduce the following notation for reflecting the arguments of functions.
\begin{notation}\label{not:phi,gamma}
For a function $f:\mR^2 \to \mH $ and a multi-index $\phi= (\phi_1, \phi_2)$ with $\phi_1, \phi_2 \in\{0,1\}$ we put
\begin{align*}
 f^{\phi} = f^{(\phi_1, \phi_2)}(x)&:=f((-1)^{\phi_{1}}x_1,(-1)^{\phi_{2}}x_2).
\end{align*}
\end{notation}
Suppose furthermore that we are dealing with two square roots of $-1$, $\mu$ and $\nu$. We introduce their anticommutator:
\begin{equation*} 
a= \{ \mu,\nu\} = \mu \nu + \nu \mu \in \mR.
\end{equation*}
This anticommutator is always in $\mR$. Indeed, when $\mu$ and $\nu$ are square roots of $-1$ they are pure quaternions (i.e. their real part is zero) and simple calculations reveal that they satisfy the following conditions:  
\begin{equation*}
\begin{aligned}
\mu &= b_1\qi + c_1\qj + d_1\qk, \qquad \quad b_1^2 + c_1^2+ d_1^2 = 1\\
\nu &= b_2\qi + c_2\qj + d_2\qk, \qquad  \quad b_2^2 + c_2^2+ d_2^2 = 1\\
a&= \{ \mu,\nu\} = \mu \nu + \nu \mu  = -2( b_1 b_2 +c_1 c_2 + d_1 d_2), 
\end{aligned}
\end{equation*}
where $\qi,\qj,\qk$ are the standard quaternion units. In particular, when $\mu$ and $\nu$ are perpendicular, then $a=0$.

We may now formulate the following theorem, see \cite{BDDS}, which expresses the Mustard convolution as a finite linear combination of classical convolutions:
\begin{equation*}
 (f* g)(x)=\int_{\mR^2}f(y) g(x-y)  dy.
\end{equation*}

\begin{theorem}
\label{th:mustclassconv}
The Mustard convolution can be expressed as, with $\{\mu,\nu\} = a \in \mR$,
\begin{align*}
f*_{\mu,\nu}g= & \frac{1}{4}\left[ f \ast g +  f \ast g^{(0,1)} + f \ast g^{(1,0)} + f \ast g^{(1,1)}   \right]\\
&+\frac{1}{4}\left[ - \nu f \ast  \nu g +   \nu f \ast \nu  g^{(0,1)} -  \nu f \ast \nu g^{(1,0)} +  \nu f \ast  \nu g^{(1,1)}    \right]\\
& +\frac{1}{4}\left[ - \mu f^{(0,1)} \ast  \mu g -   \mu f^{(0,1)}  \ast \mu g^{(0,1)} +  \mu f^{(0,1)}  \ast \mu g^{(1,0)} +  \mu f^{(0,1)}  \ast  \mu g^{(1,1)}   \right]\\
&+\frac{1}{8}\left[  - \mu \nu f^{(0,1)}  \ast  \mu \nu g +   \mu \nu f^{(0,1)}  \ast \mu \nu  g^{(0,1)} +  \mu \nu f^{(0,1)}  \ast \mu \nu g^{(1,0)} -  \mu \nu f^{(0,1)}  \ast  \mu \nu g^{(1,1)}\right]\\
&+\frac{1}{8}\left[   \nu \mu f^{(0,1)}  \ast  \mu \nu g -   \nu \mu f^{(0,1)}  \ast \mu \nu  g^{(0,1)} -  \nu \mu f^{(0,1)}  \ast \mu \nu g^{(1,0)} +  \nu \mu f^{(0,1)}  \ast  \mu \nu g^{(1,1)}\right]\\
&+\frac{a}{8}\left[   f  \ast  \mu \nu g -    f  \ast \mu \nu  g^{(0,1)} -   f \ast \mu \nu g^{(1,0)} +  f  \ast  \mu \nu g^{(1,1)}   \right.\\
&\qquad+\left.  \left(   \nu f \ast  \mu g +   \nu f \ast \mu  g^{(0,1)} -  \nu f \ast \mu g^{(1,0)} -  \nu f \ast  \mu g^{(1,1)} \right) \right.\\
&\qquad+\left.    \left( - \nu f^{(0,1)} \ast  \mu g -   \nu f^{(0,1)}  \ast \mu g^{(0,1)} +  \nu f^{(0,1)}  \ast \mu g^{(1,0)} +  \nu f^{(0,1)}  \ast  \mu g^{(1,1)} \right) \right].
\end{align*}
\end{theorem}

\begin{remark}
The formula in Theorem \ref{th:mustclassconv} simplifies dramatically for special choices of the roots $\mu$ and $\nu$. When $\mu = \pm \nu$ only four terms remain. When $\mu$ and $\nu$ are perpendicular ($a =\{\mu,\nu\}  =0$) 16 terms remain.
\end{remark}

\begin{remark}
It is possible to derive similar theorems for the case of right qFT, resp. the two-sided qFT. For the right qFT this leads to an equal number of terms. However, for the two-sided qFT the result simplifies and only 16 terms are necessary for the most general case of two roots $\mu$ and $\nu$.
\end{remark}

The main aim of the present paper is now to obtain the `inverse' result, namely to express the standard convolution $\ast$ as a finite linear combination of  Mustard convolutions $*_{\mu,\nu}$ in the following way:
\[
f \ast g = \sum_i (a_i f^{\phi_i}) *_{\mu,\nu} (b_i g^{\psi_i}), \qquad a_i, b_i \in \mH.
\]
As we will see, the proof of that statement is quite a bit more involved.

\subsection{Basic formulas}
\label{basicform}
In this section we summarize some of the basic formulas we will use to derive our results. Denote again by $\cI_{2}$ the set $\{ a \in  \mH | a^{2}=-1\}$ of geometric square roots of minus one. By $\mu,\nu$ we denote elements of  $\cI_{2}$ and by  $q,r$ elements of $ \mH$ .
\subsubsection{Commuting and anticommuting of quaternions}

Following \cite{HS, BU5}, we define the commuting ($j=0$) and the anticommuting ($j=1$) part of a quaternion $q$ with respect to $\mu$,  a square root of $-1$ as follows: 
\begin{equation*}
\begin{aligned}
q_{c^j(\mu)} &= \frac{1}{2}(q - (-1)^j \mu q \mu)  \\
q &= q_{c^0(\mu)} + q_{c^1(\mu)} = q_{-} + q_{+}  \\
q_{-} \mu &= \mu q_{-}\\
q_{+} \mu &= - \mu q_{+},
\end{aligned}
\end{equation*}
where we introduced the notations $q_{c^0(\mu)} =  q_{-} $ and  $q_{c^1(\mu)} =q_{+}$. These shorter notations are sometimes used when the quaternion we wish to commute is clear from the context.

Now suppose we are dealing with two square roots of $-1$, $\mu$ and $\nu$. Their anticommutator satisfies $a= \{ \mu,\nu\} = \mu \nu + \nu \mu \in \mR$. We can use this fact to simplify many calculations. The (anti)commuting parts of $\mu$ with respect to $\nu$ are for example:
\begin{equation*}
\begin{aligned}
\mu_{c^0(\nu)} &= \mu_{-} = \frac{1}{2}(\mu -  \nu \mu \nu ) = - \frac{a}{2} \nu\\
\mu_{c^1(\nu)} &= \mu_{+} = \frac{1}{2}(\mu +  \nu \mu \nu ) = \mu + \frac{a}{2} \nu \\
\end{aligned} 
\end{equation*}
These commuting parts have interesting properties: 
\begin{align}
\label{squaresplit}
\mu^2 = \mu_+^2 +  \mu_-^2 &= -1\\
\label{mupmprop}
\mu_+\mu_- + \mu_-\mu_+ &=0,
\end{align}
which can be derived by evaluating $- \nu$ in two different ways. Indeed, we have 
\begin{equation*}
\begin{aligned}
- \nu &=\mu^2\nu \\
&= (\mu_+ + \mu_-)^2 \nu \\
&=  (\mu_+^2 + \mu_+\mu_- + \mu_-\mu_+ +  \mu_-^2) \nu \\					
&= \nu(\mu_+^2 - \mu_+\mu_- - \mu_-\mu_+ +  \mu_-^2), 
\end{aligned}
\end{equation*}
but also
\begin{equation*}
\begin{aligned}
- \nu &=\nu \mu^2\\
&= \nu (\mu_+ + \mu_-)^2  \\
&= \nu (\mu_+^2 + \mu_+\mu_- + \mu_-\mu_+ +  \mu_-^2)\\					
\end{aligned}
\end{equation*}
which leads to \eqref{squaresplit} and \eqref{mupmprop} by equating both computations.

\subsubsection{Commuting and anticommuting of quaternions with exponentials}

In the study of the qFT we are mainly interested in commuting quaternions with exponentials such as $e^{\mu \theta} = \cos{\theta} + \mu \sin{\theta}$, with $\mu^2=-1$.  We can derive the following commutation relations for such functions:
\begin{equation}
\begin{aligned}
\label{commutingexponential}
q &= q_{-} + q_{+} = q_{c^0(\mu)} + q_{c^1(\mu)} \\
q_{\pm}e^{ \theta \mu} &= e^{\mp \theta \mu}q_{\pm}, \qquad \theta \in \mR. 
\end{aligned}
\end{equation} 
We will also frequently use its commuting and anticommuting parts: 
\begin{equation}
\label{commutingexponential2}
\begin{aligned}
e^{-\nu \theta}_{c^j(\mu)} &=  \frac{1}{2}(e^{-\nu \theta} - (-1)^j \mu (\cos{(\theta)} - \nu \sin{ (\theta)}) \mu) \\
&=  \frac{1}{2}(e^{-\nu \theta} + (-1)^j e^{\nu \theta} + a(-1)^j \mu \sin{(\theta))}\\
&=  \frac{1}{2}(e^{-\nu \theta} + (-1)^j e^{\nu \theta} + \frac{a}{2}(-1)^j e^{\mu \theta} - \frac{a}{2}(-1)^j e^{-\mu \theta})
\end{aligned}
\end{equation}
where we recall that $\{ \mu, \nu \}=a$.
For the special  case $a=0$, which we will encounter later on, the formulas reduce to: 
\begin{equation}
\label{commutingexponential1}
\begin{aligned}
e^{-\nu \theta}_{c^j(\mu)} &=  \frac{1}{2}(e^{-\nu \theta} - (-1)^j \mu e^{-\nu \theta} \mu) \\
&=  \frac{1}{2}(e^{-\nu \theta} + (-1)^j e^{\nu \theta})\\
&= \begin{cases}
                            \cos{\theta} ,&\text{ if } j =0,\\
			    -\nu \sin{\theta} ,&\text{ if } j = 1.
                           \end{cases}
\end{aligned}
\end{equation}
In this case the commuting part is the cosine, while the anticommuting part is the part with the sine.

We can apply our result with the purpose of swapping an element of $\cI_{2}$  and an exponential:  
\begin{equation}
\label{singleswap}
\begin{aligned}
\mu e^{-\nu\theta} &= (\mu_{+} + \mu_{-}) e^{-\nu\theta} \\
									 &= e^{\nu\theta} \mu_{+} + e^{-\nu \theta} \mu_{-} \\
&= e^{\nu \theta}(\mu + \frac{a}{2}\nu)  - e^{-\nu \theta}  \frac{a}{2} \nu. 
\end{aligned}
\end{equation}
Later on we will also need the swapping of two elements of $\cI_{2}$  and their exponentials: 
\begin{equation*}
\begin{aligned}
\nu \mu e^{-\mu\theta} e^{-\nu\theta} &=  [e^{\mu\theta} \nu_+ + e^{-\mu\theta} \nu_-] [e^{\nu\theta}\mu_+ + e^{-\nu\theta}\mu_-] \\
&= e^{\mu\theta} e^{\nu\theta} [ (\nu_+)_{c^0(\nu)} \mu_+ +  (\nu_+)_{c^1(\nu)} \mu_-]  +  e^{-\mu\theta} e^{\nu\theta} [(\nu_-)_{c^0(\nu)} \mu_+ + (\nu_-)_{c^1(\nu)} \mu_- ]  \\
&+  e^{\mu\theta} e^{-\nu\theta} [(\nu_+)_{c^1(\nu)} \mu_+ + (\nu_+)_{c^0(\nu)} \mu_- ] +  e^{-\mu\theta} e^{-\nu\theta} [(\nu_-)_{c^1(\nu)} \mu_+ + (\nu_-)_{c^0(\nu)} \mu_- ].
\end{aligned}
\end{equation*}
Notice that in this case the plus/minus-notation $(\mu_+)_+$ would have been ambiguous. Easy calculations yield the following: 
\begin{equation}
\label{muplusplus}
\begin{aligned}
(\nu_+)_{c^0(\nu)} &= \nu + \frac{a}{2} \mu_- \\
(\nu_+)_{c^1(\nu)} &= \frac{a}{2} \mu_+\\
(\nu_-)_{c^0(\nu)} &= - \frac{a}{2} \mu_- \\
(\nu_-)_{c^1(\nu)} &= - \frac{a}{2} \mu_+.\\
\end{aligned}
\end{equation}
Making use of \eqref{muplusplus}, and the property \eqref{squaresplit} and \eqref{mupmprop}, we can eventually simplify the expressions to obtain: 
\begin{equation}
\label{doubleswap}
\begin{aligned}
\nu \mu e^{-\mu\theta} e^{-\nu\theta} &= e^{\mu\theta} e^{\nu\theta} [ \nu \mu_+] +   e^{-\mu\theta} e^{-\nu\theta} [\frac{a}{2}]  \\
&= e^{\mu\theta} e^{\nu\theta} [ \nu \mu - \frac{a}{2} ] +   e^{-\mu\theta} e^{-\nu\theta} [\frac{a}{2}]. 
\end{aligned}
\end{equation}
\subsubsection{Properties of the quaternion Fourier transform}

We can translate the calculations on the exponentials to properties of the qFT as defined in Definition \ref{def:qft}.
A first important formula concerns the change of signs of the roots: 
\begin{equation}
\label{eq:changeofsigns}
\begin{aligned}
\cF^{(-1)^{\phi_1} \mu, (-1)^{\phi_2} \nu}(g)(u) &= \frac{1}{2\pi} \int_{\mathbb{R}^2} e^{-(-1)^{\phi_1} \mu x_1 u_1} e^{- (-1)^{\phi_2}\nu x_2 u_2} g(x) dx\\
&= \cF^{\mu, \nu}(g^{\phi} )(u),
\end{aligned}
\end{equation}
which is obtained by performing a substitution in the integral ($x \rightarrow x^{\phi}$) and using Notation \ref{not:phi,gamma}.

In the case where $\mu$ and $\nu$ are anticommuting it is also easy to obtain that passing $\mu$ and/or $\nu$ through the transform changes the sign of the corresponding exponential and thus of the arguments:  
\begin{equation}
\begin{aligned}
\label{eq:munuInsideTransform}
\mu^k\nu^l \cF^{\mu,\nu}(g)(u) &= \frac{1}{2\pi} \mu^k\nu^l   \int_{\mR^2} e^{-\mu x_1 u_1} e^{-\nu x_2 u_2} g(x) dx \\
&= \frac{1}{2\pi}  \int_{\mR^2} e^{-(-1)^l\mu x_1 u_1} e^{-(-1)^k\nu x_2 u_2} \mu^k\nu^l   g(x) dx\\ 
&= \cF^{\mu,\nu}(\mu^k\nu^l g^{(l,k)})(u).
\end{aligned}
\end{equation}

In the case where $\{\mu,\nu\} =a$, we can derive a similar property making use of equations \eqref{singleswap} and \eqref{doubleswap}. This yields after some computations: 
\begin{align}
\label{eq:munuInsideTransform2}
\begin{split}
\mu \cF^{\mu,\nu}(g)(u) &= \cF^{\mu,\nu}(\mu g^{(0,1)})(u) + \frac{a}{2}\cF^{\mu,\nu}(\nu g^{(0,1)})(u) - \frac{a}{2} \cF^{\mu,\nu}(\nu g^{(0,0)})(u) \\
\nu \cF^{\mu,\nu}(g)(u) &= \cF^{\mu,\nu}(\nu g^{(1,0)})(u) + \frac{a}{2}\cF^{\mu,\nu}(\mu g^{(1,1)})(u) - \frac{a}{2}\cF^{\mu,\nu}(\mu g^{(0,1)})(u)\\
& + \frac{a^2}{4} \cF^{\mu,\nu}(\nu g^{(1,1)})(u) - \frac{a^2}{4} \cF^{\mu,\nu}(\nu g^{(1,0)})(u)  - \frac{a^2}{4} \cF^{\mu,\nu}(\nu g^{(0,1)})(u) + \frac{a^2}{4} \cF^{\mu,\nu}(\nu g^{(0,0)})(u)
\end{split}
\end{align}
as well as
\begin{equation}
\label{eq:munuInsideTransform2_bis}
\begin{aligned}
\nu \mu \cF^{\mu,\nu}(g)(u) &=  \cF^{\mu,\nu}(\nu \mu g^{(1,1)})(u) - \frac{a}{2} \cF^{\mu,\nu}( g^{(1,1)})(u) + \frac{a}{2} \cF^{\mu,\nu}( g^{(0,0)})(u) \\
\mu \nu \cF^{\mu,\nu}(g)(u) &=  \cF^{\mu,\nu}(\mu \nu g^{(1,1)})(u) - \frac{a}{2} \cF^{\mu,\nu}( g^{(1,1)})(u) + \frac{a}{2} \cF^{\mu,\nu}( g^{(0,0)})(u).
\end{aligned}
\end{equation}

\section{Important special cases}
\label{sec:3}

Before dealing with the general case, we consider two important special cases.
\subsection{The case of equal roots}
\label{sec:eq}
When we have just one root $\mu = \nu$ the transform (\ref{eq:Transform}) reduces to: 

\begin{equation*}
\label{eq:TransformSingleRoot}
\mathcal{F}^{\mu} (f) (y) = \frac{1}{2\pi} \int_{\mathbb{R}^2} e^{-\mu (x_1 y_1 + x_2 y_2)} f(x) dx.
\end{equation*}
Now we can compute the interaction with classical convolution by applying the definitions and subsequently taking the integrals to the outside: 
\begin{equation*}
\label{startpoint1}
\begin{aligned}
\cF^{\mu}( f *g(x)) (u) &= \frac{1}{2\pi}\int_{\mR^2} e^{-\mu (x_1 u_1 +x_2 u_2)} (f *g)(x) dx \\
&=\frac{1}{2\pi} \int_{\mR^2} e^{-\mu (x_1 u_1 +x_2 u_2)} \int_{\mR^2} f(x-y)g(y) dy dx \\
&=\frac{1}{2\pi}  \int_{\mR^2} \int_{\mR^2}  e^{-\mu (x_1 u_1 +x_2 u_2)} f(x-y)g(y) dy dx \\
&= \frac{1}{2\pi} \int_{\mR^2} \int_{\mR^2}  e^{-\mu [(z_1+y_1) u_1+ (z_2+y_2) u_2 ]} f(z)g(y) dz dy \\
&= \frac{1}{2\pi} \int_{\mR^2} \int_{\mR^2}  e^{-\mu (z_1 u_1+ z_2 u_2)} e^{-\mu(y_1 u_1 +y_2 u_2 )} f(z)g(y) dz dy.
\end{aligned}
\end{equation*}
In the last step we have made the substitution $x \rightarrow z + y$. Because there is just one root, all exponential functions commute, and the exponential property holds. This will no longer be the case later on. We can continue by separating the parts in $z$ and the parts in $y$: 
\begin{equation*}
\label{commutef1}
\begin{aligned}
&\cF^{\mu}( f *g(x)) (u)\\ 
&= \frac{1}{2\pi} \int_{\mR^2} \int_{\mR^2} e^{-\mu (z_1 u_1+ z_2 u_2)}  [f(z)_{c^0(\mu)} e^{-\mu(y_1 u_1 +y_2 u_2 )} + f(z)_{c^1(\mu)} e^{\mu(y_1 u_1 +y_2 u_2 )}] g(y) dz dy\\ 
&= \frac{1}{2}\frac{1}{2\pi} \int_{\mR^2} \int_{\mR^2} e^{-\mu (z_1 u_1+ z_2 u_2)}  [ (f(z) - \mu f(z) \mu) e^{-\mu(y_1 u_1 +y_2 u_2 )} + (f(z) + \mu f(z) \mu)e^{\mu(y_1 u_1 +y_2 u_2 )}] g(y) dz dy \\
&= \frac{1}{2}\frac{1}{2\pi} \int_{\mR^2} \int_{\mR^2} \biggl( e^{-\mu (z_1 u_1+ z_2 u_2)}  f(z) e^{-\mu(y_1 u_1 +y_2 u_2 )} g(y) - e^{-\mu (z_1 u_1+ z_2 u_2)} \mu f(z) e^{-\mu(y_1 u_1 +y_2 u_2 )} \mu g(y)  \\
& \qquad + e^{-\mu (z_1 u_1+ z_2 u_2)}f(z)e^{\mu(y_1 u_1 +y_2 u_2 )} g(y)  + e^{-\mu (z_1 u_1+ z_2 u_2)}\mu f(z) e^{\mu(y_1 u_1 +y_2 u_2 )} \mu g(y)\biggr) dz dy \\
&= \frac{2 \pi}{2}  [\cF^{\mu}( f)(u)  \cF^{\mu}( g)(u) -  \cF^{\mu}( \mu f)(u)   \cF^{\mu}( \mu g)(u)  +  \cF^{\mu}( f)(u)  \cF^{\mu}( g^{(1,1)})(u)  +  \cF^{\mu}( \mu f)(u)  \cF^{\mu}( \mu g^{(1,1)})(u)].  
\end{aligned}
\end{equation*}

In a final step we apply the inverse Fourier transform to both sides. By using Definition \ref{def:conv_must} of the Mustard convolution and simplifying the notation $f *_\mu g:=f *_{\mu,\mu} g$, we rewrite the result as 
\begin{equation*}
\label{result1}
f *g= \frac{1}{2}[ f *_\mu g - \mu f *_\mu \mu g  + f *_\mu g^{(1,1)}  + \mu f *_\mu \mu  g^{(1,1)}].
\end{equation*}
This gives us our first result, which already showcases some of the important steps that will be made later on. We summarize it as follows.

\begin{theorem}
\label{th:convolutiontheorem1}
Let $f$ and $g$ be quaternion functions on $\mR^2$, and $\cF^{\mu}$ the left qFT with one root $\mu$. The classical convolution can then be expressed as a sum of Mustard convolutions by the following formula:
\begin{equation*} \label{eq:convolutionformula1}
f *g= \frac{1}{2}[ f *_\mu g - \mu f *_\mu \mu g  + f *_\mu g^{(1,1)}  + \mu f *_\mu \mu  g^{(1,1)}].
\end{equation*}
\end{theorem}

\begin{remark}
Note that one can consider in a very similar way a transform where the exponential factors are permuted in a different way, or with different signs in the exponential. 
\end{remark}

\subsection{The case of perpendicular roots}
\label{sec:perp}
Note that in this case $a=0$ and the roots $\mu$ and $\nu$ anticommute. We begin with the same steps in this more general case:
\begin{equation*}
\label{startpoint2}
\begin{aligned}
\cF^{\mu,\nu}( f *g(x)) (u) &= \frac{1}{2\pi}\int_{\mR^2} e^{-\mu x_1 u_1} e^{-\nu x_2 u_2} (f *g)(x) dx \\
&=\frac{1}{2\pi} \int_{\mR^2} e^{-\mu x_1 u_1} e^{-\nu x_2 u_2} \int_{\mR^2} f(x-y)g(y) dy dx \\
&=\frac{1}{2\pi}  \int_{\mR^2} \int_{\mR^2}  e^{-\mu x_1 u_1} e^{-\nu x_2 u_2} f(x-y)g(y) dy dx \\
&= \frac{1}{2\pi} \int_{\mR^2} \int_{\mR^2}  e^{-\mu (z_1+y_1) u_1} e^{-\nu (z_2+y_2) u_2} f(z)g(y) dz dy, 
\end{aligned}
\end{equation*}
where we again made use of the substitution $x \rightarrow z + y$. It is now our purpose to rewrite the integrand such that all exponentials with $z$ as integration variable are paired with the function $f$ and those with $y$ as integration variable with $g$. To do so, we will use the commutation relations for quaternions and exponentials as established in Section \ref{basicform}.

 First, by the use of formula \eqref{commutingexponential} we obtain for the two factors we want to interchange:
 \begin{equation}
 \label{eq:expswitch}
 e^{-\mu y_1 u_1}  e^{-\nu z_2 u_2} =  e^{-\nu z_2 u_2}_{c^0(\mu)}  e^{-\mu y_1 u_1} +  e^{-\nu z_2 u_2}_{c^1(\mu)}  e^{\mu y_1 u_1}. 
 \end{equation}
In the case of $\mu$ and $\nu$ anticommuting we can subsequently make use of the formula (\ref{commutingexponential1}). We then continue in the following way:
\begin{align*}
\cF^{\mu,\nu}( f *g(x)) (u) &= \frac{1}{2\pi} \int_{\mR^2} \int_{\mR^2}  e^{-\mu (z_1+y_1) u_1} e^{-\nu (z_2+y_2) u_2} f(z)g(y) dz dy \\
&= \frac{1}{2\pi}\int_{\mR^2} \int_{\mR^2}  \left(e^{-\mu z_1 u_1} \frac{[e^{-\nu z_2 u_2} + e^{\nu z_2 u_2}]}{2}  e^{-\mu y_1 u_1} e^{- \nu y_2 u_2} \right.\\
 &\left.+ e^{-\mu z_1 u_1} \frac{[e^{-\nu z_2 u_2} - e^{\nu z_2 u_2}]}{2}  e^{\mu y_1 u_1} e^{- \nu y_2 u_2}\right) f(z)g(y) dz dy. 
\end{align*}
Now we still have to pull the function $f$ through the exponentials. We obtain in a similar way as before:
\begin{align*}
&\cF^{\mu,\nu}( f *g(x)) (u) =  \\
& \frac{1}{2\pi} \int_{\mR^2} \int_{\mR^2} e^{-\mu z_1 u_1} \frac{[e^{-\nu z_2 u_2} + e^{\nu z_2 u_2}]}{2} \left[((f(z))_{c^0(\nu)})_{c^0(\mu)} e^{-\mu y_1 u_1} e^{- \nu y_2 u_2} + ((f(z))_{c^0(\nu)})_{c^1(\mu)} e^{\mu y_1 u_1} e^{- \nu y_2 u_2} \right.\\
& \left. + ((f(z))_{c^1(\nu)})_{c^0(\mu)} e^{-\mu y_1 u_1} e^{ \nu y_2 u_2} + ((f(z))_{c^1(\nu)})_{c^1(\mu)} e^{\mu y_1 u_1} e^{ \nu y_2 u_2} \right] g(y) dz dy\\
&+ e^{-\mu z_1 u_1} \frac{[e^{-\nu z_2 u_2} - e^{\nu z_2 u_2}]}{2}  \left[((f(z))_{c^0(\nu)})_{c^0(\mu)} e^{\mu y_1 u_1} e^{- \nu y_2 u_2} + ((f(z))_{c^0(\nu)})_{c^1(\mu)} e^{-\mu y_1 u_1} e^{- \nu y_2 u_2}  \right.\\
&+\left. ((f(z))_{c^1(\nu)})_{c^0(\mu)} e^{\mu y_1 u_1} e^{ \nu y_2 u_2} + ((f(z))_{c^1(\nu)})_{c^1(\mu)} e^{-\mu y_1 u_1} e^{+ \nu y_2 u_2} \right]g(y) dz dy.
\end{align*}
We subsequently note: 
\begin{equation*}
\begin{aligned}
\label{eq:OrderCommutingPart}
((f(z))_{c^{j_2}(\nu)})_{c^{j_1}(\mu)} &= \frac{1}{2}( (f(z))_{c^{j_2}(\nu)} - (-1)^{j_1} \mu f(z)_{c^{j_2}(\nu)} \mu) \\
&= \frac{1}{4} (f(z) - (-1)^{j_2} \nu f(z) \nu - (-1)^{j_1} \mu f(z) \mu + (-1)^{j_1 + j_2} \mu \nu f(z) \nu \mu) \\
&= \frac{1}{4} (f(z) - (-1)^{j_2} \nu f(z) \nu - (-1)^{j_1} \mu f(z) \mu + (-1)^{j_1 + j_2} \nu \mu f(z) \mu \nu )\\
&= ((f(z))_{c^{j_1}(\mu)})_{c^{j_2}(\nu)} = f(z)_{c^{[j_1,j_2]}(\mu,\nu)},
\end{aligned}
\end{equation*}
so the order in which we take (anti)commuting parts is irrelevant when $\{\mu, \nu\}=0$.

Distributing the exponentials in $\nu z_2 u_2$, which are the second factors, we get four terms. Let us name them $I_i$ with $i \in \{1,2,3,4\}$. The $I_3$ term has the same sign as $I_1$ in the $e^{-\nu z_2 u_2}$ exponential, and we therefore expect it to simplify with the $I_1$ term. Indeed, we get: 
 \begin{align*}
\begin{aligned}
I_1  &=  \frac{1}{4 \pi} \int_{\mR^2} \int_{\mR^2} e^{-\mu z_1 u_1} e^{-\nu z_2 u_2} \left[((f(z))_{c^0(\nu)})_{c^0(\mu)} e^{-\mu y_1 u_1} e^{- \nu y_2 u_2} + ((f(z))_{c^0(\nu)})_{c^1(\mu)} e^{\mu y_1 u_1} e^{- \nu y_2 u_2} \right.\\
&\left. + ((f(z))_{c^1(\nu)})_{c^0(\mu)} e^{-\mu y_1 u_1} e^{+ \nu y_2 u_2} + ((f(z))_{c^1(\nu)})_{c^1(\mu)} e^{\mu y_1 u_1} e^{ \nu y_2 u_2} \right]g(y) dz dy 
\end{aligned} \\
\begin{aligned}
I_3  &=  \frac{1}{4\pi} \int_{\mR^2} \int_{\mR^2} e^{-\mu z_1 u_1} e^{-\nu z_2 u_2} \left[((f(z))_{c^0(\nu)})_{c^0(\mu)} e^{\mu y_1 u_1} e^{- \nu y_2 u_2} + ((f(z))_{c^0(\nu)})_{c^1(\mu)} e^{-\mu y_1 u_1} e^{- \nu y_2 u_2} \right.\\
&\left.   + ((f(z))_{c^1(\nu)})_{c^0(\mu)} e^{\mu y_1 u_1} e^{ \nu y_2 u_2} + ((f(z))_{c^1(\nu)})_{c^1(\mu)} e^{-\mu y_1 u_1} e^{+ \nu y_2 u_2}\right]g(y) dz dy. 
\end{aligned}
\end{align*}
If we add these two terms and combine the exponentials in $y$ with similar signs, we notice that we have obtained sums  $((f(z))_{c^0(\nu)})_{c^0(\mu)} + ((f(z))_{c^0(\nu)})_{c^1(\mu)} = f(z)_{c^0(\nu)}$. In this way half of the terms cancel, and the following remains:
\begin{equation*}
\begin{aligned} 
I_1 + I_3 &= \frac{1}{4\pi} \int_{\mR^2} \int_{\mR^2} e^{-\mu z_1 u_1} e^{-\nu z_2 u_2} \left[f(z)_{c^0(\nu)} e^{-\mu y_1 u_1} e^{- \nu y_2 u_2} + f(z)_{c^0(\nu)} e^{\mu y_1 u_1} e^{- \nu y_2 u_2} \right. \\
&\left. + f(z)_{c^1(\nu)} e^{-\mu y_1 u_1} e^{+ \nu y_2 u_2} + f(z)_{c^1(\nu)} e^{\mu y_1 u_1} e^{ \nu y_2 u_2} \right]g(y) dz dy \\
&= \frac{\pi}{2}\left[\cF^{\mu,\nu}(f - \nu f \nu ) (u) \cF^{\mu,\nu}(g^{(0,0)})(u) + \cF^{\mu,\nu}(f  - \nu f \nu ) (u) \cF^{\mu,\nu}(g^{(1,0)})(u) \right.\\
&   \left. + \cF^{\mu,\nu}(f + \nu f \nu) (u) \cF^{\mu,\nu}(g^{(0,1)})(u)  + \cF^{\mu,\nu}(f + \nu f \nu) (u) \cF^{\mu,\nu}(g^{(1,1)})(u)\right].
\end{aligned}
\end{equation*}
Here we have identified the separate qFTs and used the change of signs property \eqref{eq:changeofsigns}. We can subsequently bring the root $\nu$ inside the qFT of $g$ by use of formula \eqref{eq:munuInsideTransform}. Let us apply the inverse transform to one of the constituting terms of $I_1+I_3$, for example $\frac{\pi}{2}\cF^{\mu,\nu}(f + \nu f \nu) (u) \cF^{\mu,\nu}(g^{(1,1)})(u)$. This gives us the following result:
\begin{equation*}
\frac{\pi}{2} {(\cF^{\mu,\nu})}^{-1}(\cF^{\mu,\nu}(f + \nu f \nu) (u) \cF^{\mu,\nu}(g^{(1,1)})(u)) = \frac{1}{4}( f \ast_{\mu,\nu} g^{(1,1)} +  \nu f \ast_{\mu,\nu} \nu g^{(0,1)}),
\end{equation*}
where we used Definition \ref{def:conv_must} of the Mustard convolution, the anticommutativity of $\mu$ and $\nu$ and formula (\ref{eq:changeofsigns}). We now conclude that applying the inverse transform to $I_1+I_3$ already gives us $8$ terms: 
\begin{equation*}
\begin{aligned}
(\cF^{\mu,\nu})^{-1} (I_1 + I_3) &= \frac{1}{4}\left[ f \ast_{\mu,\nu} g +  f \ast_{\mu,\nu} g^{(0,1)} + f \ast_{\mu,\nu} g^{(1,0)} + f \ast_{\mu,\nu} g^{(1,1)}   \right.\\
&\left.- \nu f \ast_{\mu,\nu}  \nu g +   \nu f \ast_{\mu,\nu} \nu  g^{(0,1)} -  \nu f \ast_{\mu,\nu} \nu g^{(1,0)} +  \nu f \ast_{\mu,\nu}  \nu g^{(1,1)}\right].
\end{aligned}
\end{equation*}

Next, it is easily checked that the sum of $I_2$ and $I_4$ also simplifies:
\begin{align*}
\begin{aligned}
I_2 &=  \frac{1}{4\pi} \int_{\mR^2} \int_{\mR^2} e^{-\mu z_1 u_1} e^{\nu z_2 u_2} \left[((f(z))_{c^0(\nu)})_{c^0(\mu)} e^{-\mu y_1 u_1} e^{- \nu y_2 u_2} + ((f(z))_{c^0(\nu)})_{c^1(\mu)} e^{\mu y_1 u_1} e^{- \nu y_2 u_2} \right.\\
&\left. + ((f(z))_{c^1(\nu)})_{c^0(\mu)} e^{-\mu y_1 u_1} e^{+ \nu y_2 u_2} + ((f(z))_{c^1(\nu)})_{c^1(\mu)} e^{\mu y_1 u_1} e^{ \nu y_2 u_2} \right]g(y) dz dy \\
I_4 &= - \frac{1}{4\pi} \int_{\mR^2} \int_{\mR^2} e^{-\mu z_1 u_1} e^{\nu z_2 u_2} \left[((f(z))_{c^0(\nu)})_{c^0(\mu)} e^{\mu y_1 u_1} e^{- \nu y_2 u_2} + ((f(z))_{c^0(\nu)})_{c^1(\mu)} e^{-\mu y_1 u_1} e^{- \nu y_2 u_2} \right. \\
&\left.+ ((f(z))_{c^1(\nu)})_{c^0(\mu)} e^{\mu y_1 u_1} e^{ \nu y_2 u_2} + ((f(z))_{c^1(\nu)})_{c^1(\mu)} e^{-\mu y_1 u_1} e^{+ \nu y_2 u_2}\right]g(y) dz dy.
\end{aligned}
\end{align*}
The main difference is the overall minus sign of $I_4$, which yields sums like $((f(z))_{c^0(\nu)})_{c^0(\mu)} - ((f(z))_{c^0(\nu)})_{c^1(\mu)} = - \mu (f(z))_{c^0(\nu)} \mu$. We hence obtain:
\begin{align*}
\begin{aligned}
&I_2 + I_4 \\
&= \frac{1}{4\pi} \int_{\mR^2} \int_{\mR^2} e^{-\mu z_1 u_1} e^{\nu z_2 u_2} \left[-\mu f(z)_{c^0(\nu)} \mu e^{-\mu y_1 u_1} e^{- \nu y_2 u_2} + \mu f(z)_{c^0(\nu)} \mu e^{\mu y_1 u_1} e^{- \nu y_2 u_2} \right.\\
&\left.- \mu f(z)_{c^1(\nu)} \mu e^{-\mu y_1 u_1} e^{+ \nu y_2 u_2} + \mu f(z)_{c^1(\nu)} \mu e^{\mu y_1 u_1} e^{ \nu y_2 u_2}\right]g(y) dz dy \\
&= \frac{\pi}{2} \left[\cF^{\mu,\nu}(-\mu f^{(0,1)} \mu + \mu\nu f^{(0,1)} \nu \mu) (u) \cF^{\mu,\nu}(g^{(0,0)})(u) + \cF^{\mu,\nu}(\mu f^{(0,1)} \mu   -\mu \nu f^{(0,1)} \nu \mu) (u) \cF^{\mu,\nu}(g^{(1,0)})(u) \right.\\
& \left.   + \cF^{\mu,\nu}(-\mu f^{(0,1)} \mu  -  \mu \nu f^{(0,1)}  \nu\mu) (u) \cF^{\mu,\nu}(g^{(0,1)})(u)  + \cF^{\mu,\nu}(\mu f^{(0,1)} \mu + \mu \nu f^{(0,1)} \nu \mu) (u) \cF^{\mu,\nu}(g^{(1,1)})(u)\right].
\end{aligned}
\end{align*}
Now we can again recognize the Mustard convolutions by applying the inverse transform. Using the properties \eqref{eq:changeofsigns}, \eqref{eq:munuInsideTransform} and the anticommutativity of $\mu$ and $\nu$ we finally arrive at:
\begin{equation*}
\begin{aligned}
&(\cF^{\mu,\nu})^{-1} (I_2 + I_4) \\
&= \frac{1}{4}\left[- \mu f^{(0,1)} \ast_{\mu,\nu}  \mu g -   \mu f^{(0,1)}  \ast_{\mu,\nu} \mu g^{(0,1)} +  \mu f^{(0,1)}  \ast_{\mu,\nu} \mu g^{(1,0)} +  \mu f^{(0,1)}  \ast_{\mu,\nu}  \mu g^{(1,1)}   \right.\\
&\left.  + \nu \mu f^{(0,1)}  \ast_{\mu,\nu}  \mu\nu  g -   \nu \mu f^{(0,1)}  \ast_{\mu,\nu} \mu\nu   g^{(0,1)} -  \nu \mu f^{(0,1)}  \ast_{\mu,\nu} \mu\nu  g^{(1,0)} +  \nu \mu f^{(0,1)}  \ast_{\mu,\nu}  \mu\nu  g^{(1,1)}\right].
\end{aligned}
\end{equation*}
This results in $8$ other terms. We summarize our result in the following theorem:

\begin{theorem}
\label{th:convolutiontheorem2}
Let $f$ and $g$ be quaternion functions on $\mR^2$, and $\cF^{\mu,\nu}$ the left qFT with $\mu$ and $\nu$ anticommuting roots of $-1$. The classical convolution can then be expressed as a sum of Mustard convolutions by the following formula:
\begin{equation*}
\label{eq:convolutionformula2}
\begin{aligned}
f \ast g &=  \frac{1}{4}\left[ f \ast_{\mu,\nu} g +  f \ast_{\mu,\nu} g^{(0,1)} + f \ast_{\mu,\nu} g^{(1,0)} + f \ast_{\mu,\nu} g^{(1,1)}   \right.\\
&- \nu f \ast_{\mu,\nu}  \nu g +   \nu f \ast_{\mu,\nu} \nu  g^{(0,1)} -  \nu f \ast_{\mu,\nu} \nu g^{(1,0)} +  \nu f \ast_{\mu,\nu}  \nu g^{(1,1)}  \\&- \mu f^{(0,1)} \ast_{\mu,\nu}  \mu g -   \mu f^{(0,1)}  \ast_{\mu,\nu} \mu g^{(0,1)} +  \mu f^{(0,1)}  \ast_{\mu,\nu} \mu g^{(1,0)} +  \mu f^{(0,1)}  \ast_{\mu,\nu}  \mu g^{(1,1)}   \\
& \left. + \nu \mu f^{(0,1)}  \ast_{\mu,\nu}  \mu \nu  g -   \nu \mu f^{(0,1)}  \ast_{\mu,\nu} \mu\nu   g^{(0,1)} -  \nu \mu f^{(0,1)}  \ast_{\mu,\nu} \mu\nu  g^{(1,0)} +  \nu \mu f^{(0,1)}  \ast_{\mu,\nu}  \mu\nu  g^{(1,1)}\right].
\end{aligned}
\end{equation*}
\end{theorem}

\begin{remark}
Note that in Theorem \ref{th:convolutiontheorem2} and in subsequent theorems, all terms, such as $ \nu f \ast_{\mu,\nu} \nu  g^{(0,1)}$, should be understood with parentheses $(\nu f)  \ast_{\mu,\nu} (\nu  g^{(0,1)})$. They are omitted to avoid notational overload. 
\end{remark}

\section{The general case $\{\mu,\nu\} = a$}
\label{sec:gen}
In the deduction of Theorem \ref{th:convolutiontheorem2} we assumed $\{\mu,\nu\} = a = 0$. We can now repeat the steps above while keeping $a \in \mR$ non-zero.
We need modifications in several steps. First of all the (anti)commuting part of the exponential gets extra terms. Let us recall formula \eqref{commutingexponential2}:
\begin{equation*}
\begin{aligned}
e^{-\nu z_2 u_2}_{c^j(\mu)} &=  \frac{1}{2}(e^{-\nu z_2 u_2} - (-1)^j \mu (\cos{(z_2 u_2)} - \nu \sin{ (z_2 u_2)}) \mu) \\
&=  \frac{1}{2}(e^{-\nu z_2 u_2} + (-1)^j e^{\nu z_2 u_2} + a(-1)^j \mu \sin{(z_2 u_2)})\\
&=  \frac{1}{2}(e^{-\nu z_2 u_2} + (-1)^j e^{\nu z_2 u_2} + \frac{a}{2}(-1)^j e^{\mu z_2 u_2} - \frac{a}{2}(-1)^j e^{-\mu z_2 u_2}).
\end{aligned} 
\end{equation*}
This would give terms with two exponentials in $\mu z_2 u_2$ in equation  \eqref{eq:expswitch}. Such terms are undesirable if we want to express everything in Fourier transforms $\cF^{\mu,\nu}$, but fortunately we can rewrite this expression in `good' terms (which do not depend on $a$) and `extra' terms which can be treated separately. This is done as follows: 
 \begin{equation}
 \label{eq:RewriteExpSwitch}
\begin{aligned}
 e^{-\mu y_1 u_1}  e^{-\nu z_2 u_2} &=  e^{-\nu z_2 u_2}_{c^0(\mu)}  e^{-\mu y_1 u_1} +  e^{-\nu z_2 u_2}_{c^1(\mu)}  e^{\mu y_1 u_1}\\
 &= \frac{1}{2}(e^{-\nu z_2 u_2} +  e^{\nu z_2 u_2} + \frac{a}{2} e^{\mu z_2 u_2} - \frac{a}{2} e^{-\mu z_2 u_2})e^{-\mu y_1 u_1} \\
 &+ \frac{1}{2}(e^{-\nu z_2 u_2} - e^{\nu z_2 u_2} - \frac{a}{2} e^{\mu z_2 u_2} + \frac{a}{2} e^{-\mu z_2 u_2})e^{\mu y_1 u_1} \\
 &= \text{`good terms'} + \frac{a}{4}[ e^{\mu z_2 u_2}e^{-\mu y_1 u_1} - e^{-\mu z_2 u_2}e^{-\mu y_1 u_1} - e^{\mu z_2 u_2}e^{\mu y_1 u_1} + e^{-\mu z_2 u_2} e^{\mu y_1 u_1}] \\
 &= \text{`good terms'} + \frac{a}{4}[4 \sin{(z_2 u_2)}\sin{(y_1 u_1)}]\\
 &= \text{`good terms'} + a[e^{-\nu z_2 u_2} - e^{\nu z_2 u_2}] \sin{(y_1 u_1)}  \frac{\nu}{2}\\
 &=\text{`good terms'} + \frac{a}{4} [e^{-\nu z_2 u_2} - e^{\nu z_2 u_2}][e^{-\mu  y_1 u_1} - e^{\mu y_1 u_1}] \mu \nu.
 \end{aligned}
\end{equation}

Also, when bringing the square roots inside the transform of $g$, we can no longer use formula  \eqref{eq:munuInsideTransform} but we have to use  equation \eqref{eq:munuInsideTransform2}.

Let us first focus on the modifications needed in the calculation of the `good' terms.
\subsection{Derivation of the `good' terms}

For the `good' terms we will show that bringing the roots inside the transform of $g$ does not alter the convolution formula.   
The derivation of Section \ref{sec:perp} can go unchanged up to the following point:
\begin{equation*}
\label{eq:unchanged}
\begin{aligned}
&\frac{\pi}{2}[\cF^{\mu,\nu}(f - \nu f \nu ) (u) \cF^{\mu,\nu}(g^{(0,0)})(u) + \cF^{\mu,\nu}(f  - \nu f \nu ) (u) \cF^{\mu,\nu}(g^{(1,0)})(u) \\
&    + \cF^{\mu,\nu}(f + \nu f \nu) (u) \cF^{\mu,\nu}(g^{(0,1)})(u)  + \cF^{\mu,\nu}(f + \nu f \nu) (u) \cF^{\mu,\nu}(g^{(1,1)})(u)\\
&+ \cF^{\mu,\nu}(-\mu f^{(0,1)} \mu + \mu\nu f^{(0,1)} \nu \mu) (u) \cF^{\mu,\nu}(g^{(0,0)})(u) + \cF^{\mu,\nu}(\mu f^{(0,1)} \mu   -\mu \nu f^{(0,1)} \nu \mu) (u) \cF^{\mu,\nu}(g^{(1,0)})(u) \\
&    + \cF^{\mu,\nu}(-\mu f^{(0,1)} \mu  -  \mu \nu f^{(0,1)}  \nu\mu) (u) \cF^{\mu,\nu}(g^{(0,1)})(u)  + \cF^{\mu,\nu}(\mu f^{(0,1)} \mu + \mu \nu f^{(0,1)} \nu \mu) (u) \cF^{\mu,\nu}(g^{(1,1)})(u)].
\end{aligned}
\end{equation*}
Here, we have to bring the roots inside the transform, which is now more complicated due to the lack of anticommutativity of $\mu$ and $\nu$.

To make our future arguments clearer, we will first rearrange the summands in the equation above:
\begin{equation*}
\begin{aligned}
&[\cF^{\mu,\nu}(f)\cF^{\mu,\nu}(g^{(0,0)}) + \cF^{\mu,\nu}(f)\cF^{\mu,\nu}(g^{(0,1)}) + \cF^{\mu,\nu}(f)\cF^{\mu,\nu}(g^{(1,0)}) + \cF^{\mu,\nu}(f)\cF^{\mu,\nu}(g^{(1,1)})] \\
&-  [\cF^{\mu,\nu}(\nu f \nu )  \cF^{\mu,\nu}(g^{(0,0)}) + \cF^{\mu,\nu}(\nu f \nu )  \cF^{\mu,\nu}(g^{(1,0)})] \\
&+ [ \cF^{\mu,\nu}(\nu f \nu)  \cF^{\mu,\nu}(g^{(0,1)})  +  \cF^{\mu,\nu}(\nu f \nu)  \cF^{\mu,\nu}(g^{(1,1)})]\\
& - [\cF^{\mu,\nu}(\mu f^{(0,1)} \mu)\cF^{\mu,\nu}(g^{(0,0)}) + \cF^{\mu,\nu}(\mu f^{(0,1)} \mu)\cF^{\mu,\nu}(g^{(0,1)})]\\
&+ [\cF^{\mu,\nu}(\mu f^{(0,1)} \mu)\cF^{\mu,\nu}(g^{(1,0)})+ \cF^{\mu,\nu}(\mu f^{(0,1)} \mu)\cF^{\mu,\nu}(g^{(1,1)})]\\
& + [\cF^{\mu,\nu}(\mu\nu f^{(0,1)} \nu \mu)  \cF^{\mu,\nu}(g^{(0,0)}) + \cF^{\mu,\nu}(\mu\nu f^{(0,1)} \nu \mu)  \cF^{\mu,\nu}(g^{(1,1)})]\\
& - [\cF^{\mu,\nu}(\mu\nu f^{(0,1)} \nu \mu)  \cF^{\mu,\nu}(g^{(0,1)})+ \cF^{\mu,\nu}(\mu\nu f^{(0,1)} \nu \mu)  \cF^{\mu,\nu}(g^{(1,0)})],
\end{aligned}
\end{equation*}
where we have omitted the arguments in $u$ and the factor $\pi/2$ to improve readability.

Now consider for example the terms $\cF^{\mu,\nu}(-\mu f^{(0,1)} \mu )  \cF^{\mu,\nu}(g^{(0,0)})$ and $\cF^{\mu,\nu}(-\mu f^{(0,1)} \mu )  \cF^{\mu,\nu}(g^{(0,1)})$. Using property \eqref{eq:munuInsideTransform2} we get for their sum:
\begin{equation*}
\begin{aligned}
& \cF^{\mu,\nu}(-\mu f^{(0,1)} \mu )  \cF^{\mu,\nu}(g^{(0,0)}) + \cF^{\mu,\nu}(-\mu f^{(0,1)} \mu )  \cF^{\mu,\nu}(g^{(0,1)})\\
&=-\cF^{\mu,\nu}(\mu f^{(0,1)})  \left[\cF^{\mu,\nu}(\mu g^{(0,1)}) + \frac{a}{2}\cF^{\mu,\nu}(\nu g^{(0,1)}) - \frac{a}{2} \cF^{\mu,\nu}(\nu g^{(0,0)}) \right.\\
&\left.+ \cF^{\mu,\nu}(\mu g^{(0,0)}) + \frac{a}{2}\cF^{\mu,\nu}(\nu g^{(0,0)}) - \frac{a}{2} \cF^{\mu,\nu}(\nu g^{(0,1)})\right] \\
&= -\cF^{\mu,\nu}(\mu f^{(0,1)})  \left[ \cF^{\mu,\nu}(\mu g^{0,0}) +  \cF^{\mu,\nu}(\mu g^{0,1})\right].
\end{aligned}
\end{equation*}
We observe that the terms in $a$ cancel out, and we obtain the same result as before.

This is similarly the case  for the other terms. Just to fix the idea let us consider for example the terms
$\cF^{\mu,\nu}(\mu\nu f^{(0,1)} \nu \mu)  \cF^{\mu,\nu}(g^{(0,0)})$ and $\cF^{\mu,\nu}(\mu \nu f^{(0,1)} \nu \mu)  \cF^{\mu,\nu}(g^{(1,1)})$:
\begin{equation*}
\begin{aligned}
&\cF^{\mu,\nu}(\mu\nu f^{(0,1)} \nu \mu)  \cF^{\mu,\nu}(g^{(0,0)}) + \cF^{\mu,\nu}(\mu\nu f^{(0,1)} \nu \mu)  \cF^{\mu,\nu}(g^{(1,1)})\\
&=\cF^{\mu,\nu}(\mu\nu f^{(0,1)}) \left[ \cF^{\mu,\nu}(\nu \mu g^{(1,1)}) - \frac{a}{2} \cF^{\mu,\nu}( g^{(1,1)}) + \frac{a}{2} \cF^{\mu,\nu}( g^{(0,0)}) \right.\\
&\left. +  \cF^{\mu,\nu}(\nu \mu g^{(0,0)}) - \frac{a}{2} \cF^{\mu,\nu}( g^{(0,0)}) + \frac{a}{2} \cF^{\mu,\nu}( g^{(1,1)})\right] \\ 
&= \cF^{\mu,\nu}(\mu\nu f^{(0,1)})\left[  \cF^{\mu,\nu}(\nu \mu g^{(0,0)})+   \cF^{\mu,\nu}(\nu \mu g^{(1,1)})\right].
\end{aligned}
\end{equation*}
Again the terms in $a$ cancel. We conclude that we already have the following terms in  the convolution formula: 
\begin{equation*}
\begin{aligned}
&\frac{1}{4} \left[ f \ast_{\mu,\nu} g +  f \ast_{\mu,\nu} g^{(0,1)} + f \ast_{\mu,\nu} g^{(1,0)} + f \ast_{\mu,\nu} g^{(1,1)}   \right.\\
&- \nu f \ast_{\mu,\nu}  \nu g +   \nu f \ast_{\mu,\nu} \nu  g^{(0,1)} -  \nu f \ast_{\mu,\nu} \nu g^{(1,0)} +  \nu f \ast_{\mu,\nu}  \nu g^{(1,1)}     \\
& - \mu f^{(0,1)}  \ast_{\mu,\nu}  \mu g -   \mu f^{(0,1)}  \ast_{\mu,\nu} \mu g^{(0,1)} +  \mu f^{(0,1)}  \ast_{\mu,\nu} \mu g^{(1,0)} +  \mu f^{(0,1)}  \ast_{\mu,\nu}  \mu g^{(1,1)}  \\
&\left. + \mu \nu f^{(0,1)}  \ast_{\mu,\nu}  \nu \mu g -   \mu \nu f^{(0,1)}  \ast_{\mu,\nu}\nu \mu   g^{(0,1)} -  \mu \nu f^{(0,1)}  \ast_{\mu,\nu}  \nu \mu g^{(1,0)} +  \mu \nu f^{(0,1)}  \ast_{\mu,\nu}  \nu  \mu g^{(1,1)} \right]. 
\end{aligned}
\end{equation*}
In other words, the `good' terms lead to exactly the same formula as in Theorem \ref{th:convolutiontheorem2}, except for the final swap of $\mu$ and $\nu$ in the last line.

\subsection{Derivation of the extra terms}

Next we have to continue the derivation of the extra terms in the equation \eqref{eq:RewriteExpSwitch}. Fortunately, we can use the fact that sine is real-valued and hence commutes with quaternion-valued functions. We proceed as follows: 
\begin{equation*}
\begin{aligned}
&\cF^{\mu,\nu}( f *g(x)) (u) \\
&= \frac{1}{2\pi}\int_{\mR^2} e^{-\mu x_1 u_1} e^{-\nu x_2 u_2} (f *g)(x) dx \\
&= \frac{1}{2\pi} \int_{\mR^2} \int_{\mR^2}  e^{-\mu (z_1+y_1) u_1} e^{-\nu (z_2+y_2) u_2} f(z)g(y) dz dy \\
&= \text{`good terms'} + \frac{a}{2\pi}  \int_{\mR^2} \int_{\mR^2}   e^{-\mu z_1 u_1}\sin{(z_2 u_2)} (f(z))_{c^0(\nu)} \sin{(y_1 u_1)} e^{-\nu y_2 u_2}g(y) dz dy \\
& + \frac{a}{2\pi}  \int_{\mR^2} \int_{\mR^2}   e^{-\mu z_1 u_1}\sin{(z_2 u_2)} (f(z))_{c^1(\nu)} \sin{(y_1 u_1)} e^{\nu y_2 u_2}g(y) dz dy \\
&= \text{`good terms'} + \frac{a}{8\pi}  \int_{\mR^2} \int_{\mR^2}   e^{-\mu z_1 u_1} (e^{-\nu z_2 u_2} - e^{\nu z_2 u_2}) (\nu f(z)+f(z) \nu) \sin{(y_1 u_1)} e^{-\nu y_2 u_2}g(y) dz dy \\
& + \frac{a}{8\pi}  \int_{\mR^2} \int_{\mR^2}   e^{-\mu z_1 u_1} (e^{-\nu z_2 u_2} - e^{\nu z_2 u_2}) (\nu f(z)-f(z) \nu) \sin{(y_1 u_1)} e^{\nu y_2 u_2}g(y) dz dy \\
&= \text{`good terms'} \\
&+ \frac{a}{16 \pi}  \int_{\mR^2} \int_{\mR^2} e^{-\mu z_1 u_1} (e^{-\nu z_2 u_2}-e^{\nu z_2 u_2}) \nu f(z) (e^{-\mu y_1 u_1}-e^{\mu y_1 u_1}) \mu e^{-\nu y_2 u_2} g(y) dz dy\\
& + \frac{a}{16 \pi}  \int_{\mR^2} \int_{\mR^2} e^{-\mu z_1 u_1} (e^{-\nu z_2 u_2}-e^{\nu z_2 u_2})  f(z) (e^{-\mu y_1 u_1}-e^{\mu y_1 u_1}) \mu e^{-\nu y_2 u_2} \nu g(y) dz dy\\ 
&+\frac{a}{16 \pi}  \int_{\mR^2} \int_{\mR^2} e^{-\mu z_1 u_1} (e^{-\nu z_2 u_2}-e^{\nu z_2 u_2}) \nu f(z) (e^{-\mu y_1 u_1}-e^{\mu y_1 u_1}) \mu e^{\nu y_2 u_2} g(y) dz dy\\
&- \frac{a}{16 \pi}  \int_{\mR^2} \int_{\mR^2} e^{-\mu z_1 u_1} (e^{-\nu z_2 u_2}-e^{\nu z_2 u_2})  f(z) (e^{-\mu y_1 u_1}-e^{\mu y_1 u_1}) \mu e^{\nu y_2 u_2} \nu g(y) dz dy\\
&= \text{`good terms'} + I_1 + I_2+I_3+I_4,
\end{aligned}
\end{equation*}
where we have denoted the extra terms $I_i$ with $i=1,2,3,4$.

Let us for example work out the extra terms coming from $I_1$.
We have to apply formula \eqref{singleswap} to bring $\mu$ inside the transform of $g$:
\begin{equation*}
\begin{aligned}
I_1 &= \frac{a}{16\pi}  \int_{\mR^2} \int_{\mR^2}   e^{-\mu z_1 u_1}(e^{-\nu z_2 u_2} - e^{\nu z_2 u_2}) \nu f(z) (e^{-\mu  y_1 u_1} - e^{\mu y_1 u_1}) (e^{\nu y_2 u_2}(\mu + \frac{a}{2}\nu)  - e^{-\nu y_2 u_2}  \frac{a}{2} \nu ) g(y) dz dy. 
\end{aligned}
\end{equation*}
We can now read off the extra terms coming from $I_1$ by multiplying the products out, reading the signs of the exponentials and applying the inverse qFT. This is done for the first term of $I_1$ below:
\begin{equation*}
\begin{aligned}
&{(\cF^{\mu,\nu})}^{-1} (\frac{a}{16\pi}  \int_{\mR^2} \int_{\mR^2}   e^{-\mu z_1 u_1} e^{-\nu z_2 u_2} \nu f(z) e^{-\mu  y_1 u_1} e^{\nu y_2 u_2}   \mu g(y) dz dy) \\
&=  {(\cF^{\mu,\nu})}^{-1} (\frac{a \pi}{4} \cF^{\mu,\nu}(\nu f) \cF^{\mu,-\nu}(\mu g))\\
&= \frac{a}{8}\nu f \ast_{\mu,\nu}  \mu g^{(0,1)}.
\end{aligned}
\end{equation*}
Eventually a lengthy calculation leads to the following end result for the extra terms $I_1 + I_2+I_3+I_4$:
\begin{equation*}
\begin{aligned}
 &\frac{a}{8}\left[- f \ast_{\mu,\nu}  \mu \nu g +  \nu f \ast_{\mu,\nu}  \mu g +  f \ast_{\mu,\nu}  \mu \nu g^{(0,1)} + \nu f \ast_{\mu,\nu}  \mu g^{(0,1)} \right.\\
& \qquad + f \ast_{\mu,\nu}  \mu \nu g^{(1,0)} -  \nu f \ast_{\mu,\nu}  \mu g^{(1,0)} -  f \ast_{\mu,\nu}  \mu \nu g^{(1,1)} - \nu f \ast_{\mu,\nu}  \mu g^{(1,1)}\\
& \qquad+ f^{(0,1)} \ast_{\mu,\nu}  \mu \nu g -  \nu f^{(0,1)} \ast_{\mu,\nu}  \mu g -  f^{(0,1)} \ast_{\mu,\nu}  \mu \nu g^{(0,1)} - \nu f^{(0,1)} \ast_{\mu,\nu}  \mu g^{(0,1)}\\
& \qquad\left.- f^{(0,1)} \ast_{\mu,\nu}  \mu \nu g^{(1,0)} +  \nu f^{(0,1)} \ast_{\mu,\nu}  \mu g^{(1,0)} +  f^{(0,1)} \ast_{\mu,\nu}  \mu \nu g^{(1,1)} + \nu f^{(0,1)} \ast_{\mu,\nu}  \mu g^{(1,1)}\right]\\
&+ \frac{a^2}{8}\left[f \ast_{\mu,\nu}  g -   f \ast_{\mu,\nu}  g^{(0,1)} - f \ast_{\mu,\nu}   g^{(1,0)} +  f \ast_{\mu,\nu}  g^{(1,1)} \right.\\
 &\qquad \left.- f^{(0,1)} \ast_{\mu,\nu}  g +  f^{(0,1)} \ast_{\mu,\nu}   g^{(0,1)} +  f^{(0,1)} \ast_{\mu,\nu}  g^{(1,0)} -  f^{(0,1)} \ast_{\mu,\nu}   g^{(1,1)}\right].
 \end{aligned}
\end{equation*}

We summarize the result of the derivation (sum of the `good terms' and the extra terms) in the following theorem.

\begin{theorem}
\label{th:convolutiontheorem3}
Let $f$ and $g$ be quaternion functions on $\mR^2$, and $\cF^{\mu,\nu}$ the left qFT with $\mu$ and $\nu$ roots of $-1$,  where $\{\mu, \nu \}=a$. The classical convolution can then be expressed as a sum of Mustard convolutions by the following formula:
\begin{equation*}
\label{eq:convolutiontheorem3}
\begin{aligned}
f \ast g &= \frac{1}{4} \left[ f \ast_{\mu,\nu} g +  f \ast_{\mu,\nu} g^{(0,1)} + f \ast_{\mu,\nu} g^{(1,0)} + f \ast_{\mu,\nu} g^{(1,1)}   \right. \\
&\qquad- \nu f \ast_{\mu,\nu}  \nu g +   \nu f \ast_{\mu,\nu} \nu  g^{(0,1)} -  \nu f \ast_{\mu,\nu} \nu g^{(1,0)} +  \nu f \ast_{\mu,\nu}  \nu g^{(1,1)}     \\
& \qquad- \mu f^{(0,1)}  \ast_{\mu,\nu}  \mu g -   \mu f^{(0,1)}  \ast_{\mu,\nu} \mu g^{(0,1)} +  \mu f^{(0,1)}  \ast_{\mu,\nu} \mu g^{(1,0)} +  \mu f^{(0,1)}  \ast_{\mu,\nu}  \mu g^{(1,1)}  \\
& \qquad\left.+ \mu \nu f^{(0,1)}  \ast_{\mu,\nu}  \nu \mu g -   \mu \nu f^{(0,1)}  \ast_{\mu,\nu}\nu \mu   g^{(0,1)} -  \mu \nu f^{(0,1)}  \ast_{\mu,\nu}  \nu \mu g^{(1,0)} +  \mu \nu f^{(0,1)}  \ast_{\mu,\nu}  \nu  \mu g^{(1,1)}\right] \\
 &+ \frac{a}{8}\left[- f \ast_{\mu,\nu}  \mu \nu g +  \nu f \ast_{\mu,\nu}  \mu g +  f \ast_{\mu,\nu}  \mu \nu g^{(0,1)} + \nu f \ast_{\mu,\nu}  \mu g^{(0,1)} \right. \\
&\qquad + f \ast_{\mu,\nu}  \mu \nu g^{(1,0)} -  \nu f \ast_{\mu,\nu}  \mu g^{(1,0)} -  f \ast_{\mu,\nu}  \mu \nu g^{(1,1)} - \nu f \ast_{\mu,\nu}  \mu g^{(1,1)}\\
&\qquad + f^{(0,1)} \ast_{\mu,\nu}  \mu \nu g -  \nu f^{(0,1)} \ast_{\mu,\nu}  \mu g -  f^{(0,1)} \ast_{\mu,\nu}  \mu \nu g^{(0,1)} - \nu f^{(0,1)} \ast_{\mu,\nu}  \mu g^{(0,1)}\\
&\qquad \left. - f^{(0,1)} \ast_{\mu,\nu}  \mu \nu g^{(1,0)} +  \nu f^{(0,1)} \ast_{\mu,\nu}  \mu g^{(1,0)} +  f^{(0,1)} \ast_{\mu,\nu}  \mu \nu g^{(1,1)} + \nu f^{(0,1)} \ast_{\mu,\nu}  \mu g^{(1,1)}\right]\\
&+ \frac{a^2}{8}\left[f \ast_{\mu,\nu}  g -   f \ast_{\mu,\nu}  g^{(0,1)} - f \ast_{\mu,\nu}   g^{(1,0)} +  f \ast_{\mu,\nu}  g^{(1,1)} \right.\\
 &\qquad \left.- f^{(0,1)} \ast_{\mu,\nu}  g +  f^{(0,1)} \ast_{\mu,\nu}   g^{(0,1)} +  f^{(0,1)} \ast_{\mu,\nu}  g^{(1,0)} -  f^{(0,1)} \ast_{\mu,\nu}   g^{(1,1)}\right].
\end{aligned}
\end{equation*}
\end{theorem}
In the case of $a=0$, Theorem \ref{th:convolutiontheorem3} reduces to Theorem  \ref{th:convolutiontheorem2}, which is indeed easily verified.
Note also that the derivation and the resulting expression of $f \ast g$ in terms of $ f \ast_{\mu,\nu} g$ is much more complicated than the `inverse' expression in Theorem \ref{th:mustclassconv}.

\begin{remark}
The number of terms in the right-hand side of Theorem \ref{th:convolutiontheorem3} may seem large. However, in practical applications such as image processing the formulas will often simplify. Indeed, typical masks for image filtering will satisfy certain symmetric properties. If e.g. $g^{(0,1)}=g$ the mask $g$ is an even function in the second variable, and the term in $a^2$ vanishes completely.
\end{remark}


\section{Consequences}
\label{sec:con}

A first immediate consequence is that we can now easily compute the quaternion spectrum of a classical convolution. In the case of anticommuting roots this is formulated as follows.
\begin{theorem}
\label{CCspec}
Let $f$ and $g$ be quaternion functions on $\mR^2$, and $\cF = \cF^{\mu,\nu}$ the left qFT with $\mu$ and $\nu$ anticommuting roots of $-1$. The spectrum of the classical convolution can then be expressed as follows:
\begin{equation*}
\begin{aligned}
&\cF (f \ast g )\\
 &=  \frac{\pi}{2}\left[ \cF(f) \cF( g) +  \cF(f) \cF( g^{(0,1)}) +  \cF(f) \cF( g^{(1,0)}) + \cF(f) \cF( g^{(1,1)})   \right.\\
&- \cF(\nu f)\cF(  \nu g) +   \cF(\nu f) \cF(\nu  g^{(0,1)}) - \cF( \nu f ) \cF( \nu g^{(1,0)} )+  \cF(\nu f) \cF( \nu g^{(1,1)})  \\
&- \cF(\mu f^{(0,1)}) \cF(  \mu g ) -   \cF(\mu f^{(0,1)}) \cF(\mu g^{(0,1)}) +  \cF(\mu f^{(0,1)} ) \cF( \mu g^{(1,0)}) + \cF( \mu f^{(0,1)} ) \cF( \mu g^{(1,1)} )  \\
& \left. + \cF( \nu \mu f^{(0,1)}) \cF( \mu \nu  g) -  \cF( \nu \mu f^{(0,1)} ) \cF(\mu\nu   g^{(0,1)}) -  \cF(\nu \mu f^{(0,1)}  ) \cF( \mu\nu  g^{(1,0)}) +  \cF(\nu \mu f^{(0,1)}  \cF( \mu\nu  g^{(1,1)})\right].
\end{aligned}
\end{equation*}
\end{theorem}
This result should be compared with the convolution theorem for the qFT obtained in \cite{BU5}. The advantage of our result is that it gives much better insight into the nature of the spectrum: it can be explained by suitable reflections of the original functions and multiplication by the roots $\mu$ and $\nu$. A more complete discussion of how to interpret the qFT spectrum, in the context of color image processing, is postponed to the follow-up paper \cite{SEDB}.

In data and signal processing the convolution has its counterpart in cross-correlation. It is a measure of how strongly two signals are correlated.  Classically it is defined as follows:

\begin{definition}
The classical cross-correlation of two complex functions  $f$ and $g$ defined on $\mR^2$ is given by the following integral: 
\begin{equation*}
f \star g (y) := \int_{\mR^2} f^*( x) g(x + y) dx.
\end{equation*} 
\end{definition}

The complex conjugation ($f^*$ is used to distinguish this from the quaternion conjugation $\bar{f}$) in this definition is needed to ensure a relation between the autocorrelation function $f \star f$ and the power spectrum of a signal (as required by the Wiener-Khintchine theorem). 
Notice that even in the complex case the cross-correlation $\star$ is not commutative, in contrast with the convolution $\ast$. Because of this conjugation operation there are several non-equivalent definitions of the cross correlation possible. The conjugation operation is not relevant for the derivation of our theorem and can be easily added. That is why we stick with the following definition:

\begin{definition}
The  quaternion cross-correlation of two quaternion functions on $\mR^2$, $f$ and $g$, is given by the following integral: 
\begin{equation*}
f \star g (y) := \int_{\mR^2} f( x) g(x + y) dx.
\end{equation*} 
\end{definition}
 
 With this definition we can easily observe that we can state the correlation in terms of a convolution: 
 \begin{equation*}
 f \star g  = f^{(1,1)} \ast g.
 \end{equation*}
 
 This allows to obtain results for the spectrum by using the theorems derived in this paper. As an example we will state the correlation theorem corresponding to Theorem \ref{th:convolutiontheorem2}.

\begin{theorem}
\label{th:correlationtheorem}
Let $f$ and $g$ be quaternion functions on $\mR^2$, and $\cF = \cF^{\mu,\nu}$ the left qFT with $\mu$ and $\nu$ anticommuting roots of $-1$. The spectrum of the correlation can then be expressed as a sum of products of the individual spectra by the following formula:
\begin{equation*}
\label{eq:correlationformula}
\begin{aligned}
&\cF^{\mu,\nu}(f \star g) = \\ 
& \frac{\pi}{2}\left[\cF( f^{(1,1)})  \cF( g) +  \cF(f^{(1,1)})\cF( g^{(0,1)}) + \cF(f^{(1,1)}) \cF( g^{(1,0)}) +\cF( f^{(1,1)}) \cF( g^{(1,1)})   \right.\\
&- \cF(\nu f^{(1,1)})\cF(  \nu g) +   \cF( \nu f^{(1,1)}) \cF( \nu  g^{(0,1)}) -  \cF(\nu f^{(1,1)}) \cF( \nu g^{(1,0)}) +   \cF(\nu f^{(1,1)}) \cF(  \nu g^{(1,1)})  \\&-  \cF( \mu f^{(1,0)} ) \cF(  \mu g )-   \cF(\mu f^{(1,0)}  ) \cF(  \mu g^{(0,1)}) +  \cF(\mu f^{(1,0)}  ) \cF(  \mu g^{(1,0)}) +  \cF(\mu f^{(1,0)}  ) \cF(   \mu g^{(1,1)} )  \\
&\left.  +\cF( \nu \mu f^{(1,0)}  ) \cF(   \mu \nu  g) -   \cF(\nu \mu f^{(1,0)}  ) \cF(  \mu\nu   g^{(0,1)}) -  \cF(\nu \mu f^{(1,0)}  ) \cF(  \mu\nu  g^{(1,0)}) + \cF( \nu \mu f^{(1,0)}  ) \cF(   \mu\nu  g^{(1,1)})\right].
\end{aligned}
\end{equation*}
\end{theorem}

\section{Conclusions and outlook}
\label{sec:6}

In this paper we have established formulas to express the classical convolution product of two quaternion functions as a linear combination of Mustard convolutions for the left qFT. The computations leading to this result can be adapted without additional difficulty to other versions of the qFT. In a subsequent paper we plan to show how to apply our results to the design of quaternion filters for color images. We also plan to extend our results to higher dimension, where the role of the qFT is taken over by a hypercomplex Fourier transform defined using an arbitrary number of roots of -1.

\section*{Acknowledgements}

The work of N. De Schepper was supported in part by the Fund for Scientific Research-Flanders (FWO-V), grant G.0116.13N.  A visit of T.A. Ell and S.J. Sangwine to Ghent University was funded through UGent BOF starting grant 01N01513.


\end{document}